\documentclass[12pt, a4paper]{article}

\usepackage[utf8]{inputenc}
\usepackage{titlesec}
\usepackage[english]{babel}
\usepackage{graphicx}
\usepackage{float}
\usepackage{fancyhdr}
\usepackage[matrix,arrow,curve]{xy}
\usepackage{pstricks} 
\usepackage{amsmath,amsfonts,verbatim,afterpage,theorem,euscript,mathrsfs,amssymb}
\usepackage{amsfonts}
\usepackage{amssymb}
\usepackage{array}
\usepackage{dsfont}
\usepackage{hyperref}
\usepackage{authblk}

\newtheorem{Definition}{Definition}[section]
\newtheorem{Proposition}{Proposition}[section]

\newtheorem{Lemma}{Lemma}[section]
\newtheorem{Theorem}{Theorem}

\def \vu{\textbf{u}}

\def \vw{\textbf{w}}

\newcommand\vN{{\bf  \nabla}}

\def \F{\mathbb{F}}

\def \Rd{\mathbb{R}^{d}}
\def \Rt{\mathbb{R}^{3}}
\def \Rdo{\mathbb{R}^{2}}

\def \pv{{\bf{Proof.}}~}

\newcommand \Endproof{\hfill $\diamond$}

\begin{document}

\title{Characterisation of the pressure term in the incompressible Navier--Stokes equations on the whole space}
\author{Pedro Gabriel Fern\'andez-Dalgo\footnote{LaMME, Univ Evry, CNRS, Universit\'e Paris-Saclay, 91025, Evry, France } \footnote{e-mail : pedro.fernandez@univ-evry.fr} ,  Pierre Gilles Lemari\'e--Rieusset\footnote{LaMME, Univ Evry, CNRS, Universit\'e Paris-Saclay, 91025, Evry, France} \footnote{e-mail : pierregilles.lemarierieusset@univ-evry.fr}}
\date{}\maketitle

\maketitle

\begin{abstract}
We characterise the pressure term in the incompressible 2D and 3D Navier--Stokes equations for solutions defined on the whole space.
 \end{abstract}
 

\section*{Introduction}

In the context of the Cauchy initial value problem for Navier--Stokes equations on $\mathbb{R}^d$  (with $d=2$ or $d=3$)
\begin{equation*}   \left\{ \begin{matrix} \partial_t \vu= \Delta \vu  -(\vu\cdot \vN)\vu- \vN p +\vN\cdot \mathbb{F} \cr \cr \vN\cdot \vu=0,    \phantom{space space} \vu(0,.)=\vu_0
 \end{matrix}\right.\end{equation*}
an important problem is to propose  a   formula  for the gradient of the pressure, which is an auxiliary unknown  (usually interpreted as a Lagrange multiplier for the constraint of incompressibility).

As we shall not assume differentiability of $\vu$ in our computations, it is better to write the equations as
\begin{equation*}  (NS) \left\{ \begin{matrix} \partial_t \vu= \Delta \vu  -\vN\cdot (\vu \otimes\vu)- \vN p +\vN\cdot \mathbb{F} \cr \cr \vN\cdot \vu=0,    \phantom{space space} \vu(0,.)=\vu_0
 \end{matrix}\right.\end{equation*}

Taking the Laplacian of equations (NS), since we have for a vector field $\vw$ the identity
$$ -\Delta \vw= \vN\wedge(\vN\wedge \vw)- \vN(\vN\cdot\vw) $$ 
we get the equations
$$ \partial_t\Delta u= \Delta^2 u+\vN\wedge(\vN\wedge (\vN\cdot (\vu\otimes\vu-\mathbb{F})))$$ and $$ 0=-\Delta \vN p- \vN(\vN\cdot (\vN\cdot (\vu\otimes\vu-\mathbb{F}))=-\Delta \vN p-\vN( \sum_{1\leq i,j\leq d} \partial_i\partial_j(u_iu_j-F_{i,j})).$$
Thus, the rotational-free unknown $\vN p$ obeys a Poisson equation. If $G_d$ is the fundamental solution of the operator $-\Delta$ :
$$ G_2=\frac 1{2\pi}\ln(\frac 1{\vert x\vert}),\quad G_3=\frac 1{4\pi\vert x\vert}$$ (which satisifies $-\Delta G_d=\delta$), we \textit{ formally } have
\begin{equation} \label{pressure}\vN p= G_d*\vN( \sum_{1\leq i,j\leq d} \partial_i\partial_j(u_iu_j-F_{i,j}))+ H \end{equation} with $\Delta H=0$. In the litterature, one usually finds the assumption that $\vN p$ vanishes at infinity and this is read as $H=0$. Equivalently, this is read as $$ \partial_t\vu =G_d*\vN\wedge(\vN\wedge \partial_t \vu);$$ the operator
$$ \mathbb{P}=G_d*\vN\wedge(\vN\wedge .)$$ is called the Leray projection operator and the decomposition (\textit{when justified})
$$ \vw = \mathbb{P}\vw +  G_d*\vN(\vN\cdot\vw)$$ the Hodge decomposition of the vector field $\vw$.

Hence, an important issue when dealing with the Navier--Stokes equations is to study whether in formula (\ref{pressure}) the first half of the right-hand term is well-defined, and if so which values  the second half (the harmonic part $H$) may have.

In order to give some meaning to the formal convolution $G_d*\vN\partial_i\partial_j (u_i u_j)$ or to $\left(\vN\partial_i\partial_j G_d\right)*(u_iu_j)$, we should require $u_i$ to be locally $L^2_t L^2_x$ (in order to define $u_iu_j$ as a distribution) and to have small increase at infinity, since the distribution $\vN\partial_i\partial_j G$ has small decay at infinity (it belongs  to $L^1\cap L^\infty$ far from the origin and is $O(\vert x\vert^{ -(d+1)})$).
Thus, we will focus on solutions $\vu$ that belong to $L^2((0,T), L^2(\mathbb{R}^d,w_{d+1}\, dx))$ where
$$ w_\gamma(x)=(1+\vert x\vert)^{-\gamma}.$$
 We shall recall various examples  (from recent or older litterature) of solutions belonging to the space  $L^2((0,T), L^2(\mathbb{R}^d, w_\gamma\, dx))$ (with $\gamma\in \{d,d+1\}$). As $F_{i,j}$ plays a role similar to $u_iu_j$, we shall assume that $\mathbb{F}\in L^1((0,T),L^1(\mathbb{R}^d, w_\gamma\, dx))$.

\section{Main results.}

First, we precise the meaning of $\vN p$ in equations (NS) :

\begin{Lemma}\label{curlfree}
  Consider the dimension $d \in \{ 2 , 3\}$ and $\gamma \geq 0$. Let $0<T<+\infty$. Let   $\mathbb{F}$ be  a tensor $\F(t,x)=\left(F_{i,j}(t,x)\right)_{1\leq i,j\leq d}$ such that $\F$ belongs to $L^1((0,T),L^1(\Rd, w_\gamma\, dx))$, and let $\vu$ be a vector field $\vu(t,x)=(u_i(t,x))_{1\leq i\leq d}$ such that  $\vu$ belongs to $L^2 ((0,T),L^2(\Rd, w_\gamma\, dx))$ and $\vN\cdot\vu=0$. Define the distribution $S$ by
  $$ {\bf S}=  \Delta \vu -\vN\cdot(\vu\otimes\vu-\mathbb{F})-\partial_t\vu.$$ Then the following assertions are equivalent :
  \\ (A) $S$ is curl-free : $\vN\wedge {\bf S}=0$.
  \\ (B) There exists a distribution $p\in \mathcal{D}'((0,T)\times\mathbb{R}^d)$ such that $S=\vN p$.
\end{Lemma}
  \begin{Theorem}
  \label{th}
  Consider the dimension $d \in \{ 2 , 3\}$. Let $0<T<+\infty$. Let   $\mathbb{F}$ be  a tensor $\F(t,x)=\left(F_{i,j}(t,x)\right)_{1\leq i,j\leq d}$ such that $\F$ belongs to $L^1((0,T),L^1(\Rd, w_{d+1}\, dx))$.
  Let $\vu$ be a solution of the following problem 
 \begin{equation}  \label{advection difusion free} \left\{ \begin{matrix} \partial_t \vu= \Delta \vu  - \vN \cdot  (\vu\otimes\vu)- {\bf S}+\vN\cdot \mathbb{F} \cr \cr \vN\cdot \vu=0, \quad \vN\wedge {\bf S}=0,    \phantom{space space} 
 \end{matrix}\right.\end{equation}
such that :
  $\vu$ belongs to $L^2 ((0,T),L^2_{w_{d+1}}(\Rd))$, and ${\bf S}$ belongs to $\mathcal{D}'( (0,T) \times \mathbb{R} ^d )$. \\
  
  Let us choose   $\varphi \in \mathcal{D}(\Rd)$   such that $\varphi(x )=1$ on a neighborhood of $0$ and define  $$A_{i,j,\varphi}=(1-\varphi) \partial_i \partial_j G_d. $$ Then,   there exists $g(t)\in L^1((0,T))$ such that
   $${\bf S}=\vN p_\varphi +\partial_tg$$ with \begin{align*} p_\varphi= &  \sum_{i,j}(\varphi \partial_i \partial_j G_d) * (u_iu_j-F_{i,j}) \cr &+   \sum_{i,j}   \int (A_{i,j,\varphi}(x-y)-A_{i,j,\varphi}(-y))  (u_i(t,y)u_j(t,y)-F_{i,j}(t,y))\, dy. \end{align*}
   
   Moreover, 
   \begin{itemize}
   \item[$\bullet$] $\vN p_\varphi$ does  not depend on the choice of $\varphi$ : if we change $\varphi$ in $\psi$, then
   $$ p_\varphi(t,x)-p_\psi(t,x)=  \sum_{i,j}   \int (A_{i,j,\psi}(-y)-A_{i,j,\varphi}(-y))  (u_i(t,y)u_j(t,y)-F_{i,j}(t,y))\, dy. $$
   \item[$\bullet$] $\vN p_\varphi$ is the unique solution of the Poisson problem
   $$ \Delta \vw= - \vN(\vN\cdot (\vN\cdot (\vu\otimes\vu-\mathbb{F})) $$
   with 
   $$ \lim_{\tau\rightarrow +\infty} e^{\tau\Delta}\vw =0 \text{ in }\mathcal{D}'.$$
   \item if $\mathbb{F}$ belongs more precisely to  $L^1 ((0,T),L^1_{w_{d}}(\Rd))$
and $\vu$ belongs to $L^2 ((0,T),L^2_{w_{d}}(\Rd))$, then $g=0$ and $\vN p_\varphi=\vN p_0$ where
  \begin{equation*}
    p_0=   \sum_{i,j}(\varphi \partial_i \partial_j G_d) * (u_iu_j-F_{i,j}) +   \sum_{i,j}((1-\varphi) \partial_i \partial_j G_d) * (u_iu_j-F_{i,j}).
  \end{equation*} 
  ($p_0$ does not actually depend on $\varphi$ and could have been defined as 
  $p_0=   \sum_{i,j}(  \partial_i \partial_j G_d) * (u_iu_j-F_{i,j})$.) 
   \end{itemize}
      \end{Theorem}

  When $\mathbb{F}=0$, the   case $g\neq 0$    can easily be reduced to a change of referential, due to the extended Galilean invariance of the Navier--Stokes equations~:

  \begin{Theorem}
  \label{thbis} 
  Consider the dimension $d \in \{ 2 , 3\}$. Let $0<T<+\infty$.    Let $\vu$ be a solution of the following problem 
 \begin{equation}  \label{advection difusion free} \left\{ \begin{matrix} \partial_t \vu= \Delta \vu  - \vN \cdot  (\vu\otimes\vu)- {\bf S}  \cr \cr \vN\cdot \vu=0, \quad \vN\wedge {\bf S}=0,  \quad    \vu(0,x)=\vu_0(x)
 \end{matrix}\right.\end{equation}
such that :
  $\vu$ belongs to $L^2 ((0,T),L^2_{w_{d+1}}(\Rd))$, and ${\bf S}$ belongs to $\mathcal{D}'( (0,T) \times \mathbb{R} ^d )$. \\
  
  Let us choose   $\varphi \in \mathcal{D}(\Rd)$   such that $\varphi(x )=1$ on a neighborhood of $0$ and define  $$A_{i,j,\varphi}=(1-\varphi) \partial_i \partial_j G_d. $$ We decompose  ${\bf S}$ into 
   $${\bf S}=\vN p_\varphi +\partial_tg$$ with \begin{align*} p_\varphi= &  \sum_{i,j}(\varphi \partial_i \partial_j G_d) * (u_iu_j) \cr &+   \sum_{i,j}   \int (A_{i,j,\varphi}(x-y)-A_{i,j,\varphi}(-y))  (u_i(t,y)u_j(t,y) )\, dy \end{align*} and $$g(t)\in L^1((0,T)).$$\\
   
   Let us define $$E(t)= \int_0^t g (\lambda) d\lambda$$ and $$\vw(t,x ) =\vu(t, x -E(t)) + g(t).$$ Then, $\vw$ is a solution of  the Navier--Stokes problem
    \begin{equation}   \left\{ \begin{matrix} \partial_t \vw= \Delta \vw  - \vN \cdot  (\vw\otimes\vw)- \vN q_\varphi \cr \cr \vN\cdot \vw=0, \quad \quad    \vw(0,x)=\vu_0(x) \cr\cr  q_\varphi=   {\displaystyle \sum_{i,j} } (\varphi \partial_i \partial_j G_d) * (w_iw_j)  +    {\displaystyle \sum_{i,j} }    \int (A_{i,j,\varphi}(x-y)-A_{i,j,\varphi}(-y))  (w_i(t,y)w_j(t,y) )\, dy
 \end{matrix}\right.\end{equation}  
  \end{Theorem}

\section{Curl-free vector fields.}
In this section we prove Lemma \ref{curlfree} with simple arguments  :\\

\pv{} We take a partition of unity on $(0,T)$ $$\sum_{j\in\mathbb{Z}} \omega_j=1 $$ with $\omega_j$ supported in $(2^{j-2}T, 2^jT)$ for $j<0$, in $(T/4, 3T/4)$ for $j=0$ and in   $(T-2^{-j}T, T-2^{-(j+2)}T)$ for $j>1$. We define
$$ V_j=-\omega_j \vu+\int_0^t \omega_j \Delta \vu -  \omega_j \vN\cdot(\vu\otimes\vu-\mathbb{F})+(\partial_t\omega_j)\vu\, ds.$$
Then $V_j$  is a sum of the form $A+ \Delta B+\vN\cdot C+D$ with $A$, $B$, $C$ and $D$ in $L^1((0,T),L^1(\Rd, w_{d+1}\, dx))$; thus, by Fubini's theorem, we may see it as a time-dependent tempered distribution.  Moreover, $\partial_t V_j=\omega_j {\bf S}$, $V_j$  is equal to $0$ for $t$ in a neighbourhood of $0$, and $\vN\wedge V_j=0$. Moreover, ${\bf S}=\sum_{j\in\mathbb{Z}} \partial_t V_j$.

We choose $\Phi\in\mathcal{S}(\Rd)$ such that the Fourier transform of $\Phi$ is compactly supported and is equal to $1$ in the neighbourhood of $0$. Then $\Phi*V_j$ is well-defined  and $\vN\wedge(\Phi* V_j)=0$. We define
$$ X_j =\Phi*V_j \text{ and } Y_j=V_j-X_j.$$ We have $$Y_j=\vN \left( \frac 1\Delta \vN\cdot Y_j\right)$$ and (due to Poincar\'e's lemma)
$$ X_j=\vN(\int_0^1 x\cdot X_j(t,\lambda x)d\lambda)$$
We find ${\bf S}=\vN p$ with
$$ p=\partial_t\sum_{j\in\mathbb{Z}}( \int_0^1 x\cdot X_j(t,\lambda x)d\lambda+ \frac 1\Delta \vN\cdot Y_j).$$
\Endproof

\section{ The Poisson problem $-\Delta U = \partial_k  \partial_i \partial_j  h $}

We first consider a simple Poisson problem :

\begin{Proposition}
\label{poisson}
Let $h \in L^1 (\Rd, (1+|x|)^{-(d+1)} dx)$ then 
\begin{equation*}
    U = U_1 + U_2 = (\partial_k (\varphi \partial_i\partial_j G_d))*h + \partial_k( (1- \varphi) \partial_i \partial_j G_d)* h.
\end{equation*}
is a distribution such that $U_2$ belongs to  $ L^1 (\Rd, (1+|x|)^{-(d+1)} dx)$ and $U$ is a solution of the problem
\begin{equation}
    -\Delta U = \partial_k \partial_i \partial_j  h.
\end{equation} More precisely,  $U$ is the unique solution in $\mathcal{S}'$ such that $\lim_{\tau\rightarrow 0} e^{\tau\Delta}U=0$ in $\mathcal{S}'$.
\end{Proposition}

\pv 
We may write  $\partial_j G_d$ as
\begin{equation*}
    \partial_j G_d =- \int_0^{+\infty} \partial_j W_t dt
\end{equation*}
where $W_t (x)$ is the heat kernel $W_t (x) = ({4 \pi} t)^{-\frac{d}{2}} e^{-\frac{|x|^2}{4t}}$,
so that on $\Rd \setminus \{0\}$, we have
\begin{equation*}
    \partial_j G_d = c_d \frac{x_i}{|x|^d} \text{ with } c_d= \frac 1{2(4\pi)^{d/2}} \int_0^{+\infty} e^{-\frac{1}{4u}} \frac{du}{u^{\frac{d+2}{2}}}
\end{equation*}

The first part defining $U$, $U_1=(\partial_k (\varphi \partial_i\partial_j G_d))*h $,  is well defined,  since $\partial_k (\varphi \partial_i\partial_j G_d)$ is a  compactly supported distribution. To control $U_2$, we  write
\begin{align*}
    \int & \int \frac{1}{(1+|x|)^{d+1}} |\partial_k (1- \varphi) \partial_i \partial_j G_d(x-y)| |h(y)| dy dx \\
    & \leq \int \int \frac{1}{(1+|x|)^{d+1}} \frac{C}{(1+|x-y|)^{d+1}} |h(y)| dy dx \\
    & \leq C'  \int \frac{1}{(1+|y|)^{d+1}} |h(y)| dy  \\
\end{align*}
since
\begin{align*}
    &\int \frac{1}{(1+|x|)^{d+1}} \frac{1}{(1+|x-y|)^{d+1}}  dx  \\ 
    & \leq \int_{|x|>\frac{|y|}{2}} \frac{1}{(1+|x|)^{d+1}} \frac{1}{(1+|x-y|)^{d+1}} dx    + \int_{|x-y|>\frac{|y|}{2}} \frac{1}{(1+|x|)^{d+1}}  \frac{1}{(1+|x-y|)^{d+1}} dx \\
    & \leq \frac{2^{d+1}}{(1+|y|)^{d+1}} \int \frac{1}{(1+|x-y|)^{d+1}} dx    + \frac {2^{d+1}}{(1+|y|)^{d+1}}  \int \frac{1}{(1+|x|)^{d+1}} dx \\
    & \leq C \frac{1}{(1+|y|)^{d+1}}.
\end{align*}

Now that we know that $U$ is well defined, we may compute $-\Delta U$. $-\Delta U_1$ is equal to $$(-\Delta\partial_k (\varphi \partial_i\partial_j G_d))*h= \partial_k (\varphi\partial_i\partial_jh) -\partial_k((\Delta\varphi)\partial_i\partial_jG_d)*h-2\sum_{1\leq l\leq d} \partial_k((\partial_l\varphi)\partial_l\partial_i\partial_jG_d)*h.$$
For computing $-\Delta U_2$, we see that we can differentiate under the integration sign and find
$$-\Delta U_2= = \partial_k ((1-\varphi)\partial_i\partial_jh) + \partial_k((\Delta\varphi)\partial_i\partial_jG_d)*h+2\sum_{1\leq l\leq d} \partial_k((\partial_l\varphi)\partial_l\partial_i\partial_jG_d)*h.$$
Thus, $U$ is a solution of the Poisson problem.\\

Computing $e^{\tau\Delta}U$, we find that
$$e^{\tau\Delta}U =(e^{\tau\Delta} \partial_k\partial_i\partial_jG_d)* h$$ and thus
$$ \vert e^{\tau\Delta}U (x)\vert\leq C \int \frac 1{(\sqrt \tau+ \vert x-y\vert)^{d+1}} \vert h(y)\vert\, dy.$$ By the dominated convergence theorem, we get that $\lim_{\tau\rightarrow 0} e^{\tau\Delta}U=0$ in     $ L^1 (\Rd, (1+|x|)^{-(d+1)} dx)$. If $V$ is another solution of the same Poisson problem with $V\in\mathcal{S}'$ and  $\lim_{\tau\rightarrow 0} e^{\tau\Delta}V=0$ in   $\mathcal{S}'$, then $\Delta(U-V)=0$ and $U-V\in\mathcal{S}'$, so that $U-V$ is a polynomial; with the assumption that $\lim_{\tau\rightarrow 0} e^{\tau\Delta}(U-V)=0$, we find that this polynomial is equal to $0$.\Endproof\\

If we have better integrability on $h$, then of course we have better integrability of $U_2$. For instance, we have : 
\begin{Proposition}
\label{poissonbis}
Let $h \in L^1 (\Rd, (1+|x|)^{-d} dx))$ then 
\begin{equation*}
   U_2 =  \partial_k( (1- \varphi) \partial_i \partial_j G_d)* h.
\end{equation*}
  belongs to  $  L^1 (\Rd, (1+|x|)^{-d}))$.\end{Proposition}

\pv{} 
 
We write
\begin{align*}
    \int & \int \frac{1}{(1+|x|)^{d}} |\partial_k( (1- \varphi) \partial_i \partial_j G_d(x-y))| |h(y)| \, dy \,  dx \\
    & \leq C\int \int \frac{1}{(1+|x|)^{d}} \frac{1}{(1+|x-y|)^{d+1}} |h(y)| \, dy \, dx.\end{align*}
For $|y|<1$, we have
\begin{equation*}
\label{c1p1}
    \int \frac{1}{(1+|x|)^{d+1}} \frac{1}{(1+|x-y|)^{d+1}}  dx  \leq     \int \frac{1}{(1+|x|)^{d+1}}    dx  \leq  C
\end{equation*}
and for $|y|>1$, as the real number $\int_{|x|<\frac{1}{2}} \frac{1}{|x|^{d-1}}
    \frac{1}{|x-\frac{y}{|y|}|^{2}} dx$ is finite and does not depend on $y$, we can write
\begin{align*}
    &\int \frac{1}{(1+|x|)^{d}} \frac{1}{(1+|x-y|)^{d+1}}  dx  \\ 
    & \leq \int_{|x|>\frac{|y|}{2}} \frac{1}{(1+|x|)^{d}} \frac{1}{(1+|x-y|)^{d+1}} dx    + \int_{|x|<\frac{|y|}{2}} \frac{1}{(1+|x|)^{d}}  \frac{1}{(1+|x-y|)^{d+1}} dx \\
    & \leq \frac{2^d}{(1+|y|)^{d}} \int \frac{1}{(1+|x-y|)^{d+1}} dx    + \frac{2^{d-1}}{(1+|y|)^{d-1}}  \int_{|x|<\frac{|y|}{2}} \frac{1}{|x|^{d-1}}
    \frac{1}{|x-y|^{2}} dx \\
    & \leq C \frac{1}{(1+|y|)^{d}} + C \frac{1}{(1+|y|)^{d-1}}\frac{1}{|y|}  \int_{|x|<\frac{1}{2}} \frac{1}{|x|^{d-1}}
    \frac{1}{|x-\frac{y}{|y|}|^{2}} dx \\
    & \leq C' \frac{1}{(1+|y|)^{d}}.
\end{align*} This concludes the proof.
\Endproof{}

\section{ The Poisson problem $-\Delta  V =  \partial_i \partial_j  h $}

\begin{Proposition}
\label{poisson2}
Let $h \in L^1 ((1+|x|)^{-d-1} dx)$ and $A_\varphi = (1-\varphi) \partial_i \partial_j G_d $ then 
\begin{align*}
    V = V_1 + V_2 =  (\varphi \partial_i \partial_j G_d) * h + \int ( A_\varphi(x-y) - A_\varphi(-y) )  h(y) dy
\end{align*}
is a distribution such that $V_2$ belongs to  $  L^1 ((1+|x|)^{-\gamma})$, for $\gamma > d+1$, and $V$ is a solution of the problem
\begin{equation}
    -\Delta   V =  \partial_i \partial_j h.
\end{equation}

\end{Proposition}

\pv{}
We know that $V_1$ is well defined since $ \varphi \partial_i \partial_j G_d $ is a supported compactly distribution, and we will verify that $V_2$ is well defined. 

We have
\begin{align*}
    \int & \int \frac{1}{(1+|x|)^{\gamma}} |A_\varphi(x-y) - A_\varphi (-y)|  dx \\
    & \leq C \| A_\varphi \|_{L^{\infty}}  \int \frac{1}{(1+|x|)^{\gamma}} dx  \\
\end{align*}

For $|y|>1$, we have by the mean value inequality
\begin{align*}
    \int_{|x|<{\frac{|y|}{2}}}  \frac{1}{(1+|x|)^{\gamma}} |A_\varphi(x-y) - A_\varphi (-y)|  dx \leq  C \frac{1}{|y|^{d+1}}  \int_{|x|<{\frac{y}{2}}}  \frac{|x|}{(1+|x|)^{\gamma}}  dx
\end{align*}
and we can control the other part as follows
\begin{align*}
    \int_{|x|>{\frac{|y|}{2}}}  \frac{1}{(1+|x|)^{\gamma}} |A_\varphi(-y) |  dx 
    \leq  C \frac{1}{|y|^{d}}  \int_{|x|>{\frac{y}{2}}}  \frac{1}{(|x|)^{\gamma}}  \leq  C \frac{1}{|y|^{\gamma}}
\end{align*}
and for $\varepsilon> 0$ such that $\gamma - \varepsilon \geq d+1$, we have
\begin{align*}
    \int_{|x|>{\frac{|y|}{2}}}  \frac{1}{(1+|x|)^{\gamma}} |A_\varphi(x-y) |  dx 
    &\leq   C  \int_{|x|>{\frac{|y|}{2}}}  \frac{1}{|x|^{\gamma}} \frac{1}{(1+|x-y|)^{d}}  dx\\
    &\leq   C  \int_{|x|>{\frac{|y|}{2}}}  \frac{1}{|x|^{\gamma}} \frac{1}{|x-y|^{d-\varepsilon}}  dx \\
    &\leq   C \frac{1}{|y|^{d+1}}  \int_{|x|> \frac{1}{2}} \frac{1}{|x|^{\gamma}}
    \frac{1}{|x-\frac{y}{|y|}|^{d-\varepsilon}} dx  \\ &\leq C' \frac 1{\vert y\vert^d}.
\end{align*}
\Endproof{}

  \section{Proof of Theorems \ref{th} and \ref{thbis}.}

We may now prove Theorem \ref{th} :\\

\pv{} 
Taking the divergence of
  \begin{equation*}
      \partial_t  \vu = \Delta \vu - \vN \cdot (\vu\otimes\vu) -{\bf S} + \vN \cdot \F,
  \end{equation*}
 we obtain
  \begin{equation*}
      -\sum_{i,j} \partial_i \partial_j (u_i u_j) +\sum_{i,j} \partial_i \partial_j F_{i,j} - \vN \cdot {\bf S}  =0
  \end{equation*}
  and 
  \begin{equation*}
      - \Delta {\bf S}  = \vN (\sum_{i,j} \partial_i \partial_j (u_i u_j- F_{i,j} )  ). 
  \end{equation*}
  
We write $h_{i,j} = u_i u_j - F_{i,j}$, and    
  $A_\varphi = (1-\varphi) \partial_i \partial_j G_d $. By Proposition \ref{poisson2}, we can define
\begin{align*}
    p_\varphi=  \sum_{i,j}(\varphi \partial_i \partial_j G_d) * h_{i,j} + \sum_{i,j} \int ( A_{i,j,\varphi}(x-y) - A_{i,j\varphi}(-y) )  h_{i,j}(y) dy
\end{align*}
and
\begin{equation*}
    U= U_1+ U_2 = \vN \sum_{i,j}(\varphi \partial_i \partial_j G_d) * h_{i,j} + \vN \sum_{i,j}((1-\varphi) \partial_i \partial_j G_d) * h_{i,j} = \vN p_\varphi.
\end{equation*}
Let $\tilde U= {\bf S} -U$.   First, we remark that $\Delta  U = \Delta {\bf S} $ so that $\Delta \tilde U = 0$, hence $\tilde U $ is   harmonic in the space variable.

On the other hand, for a test function $\alpha \in \mathcal{D}(\mathbb{R})$ such that $\alpha (t)= 0$  for all $|t| \geq \varepsilon$, and a test function $\beta\in\mathcal{D}(\mathbb{R}^3)$,  and  for $t \in (\varepsilon, T-\varepsilon)$, we have 
 \begin{equation*} \begin{split} \tilde U(t) *_{t,x} (\alpha\otimes\beta)  =& ( \vu *(-\partial_t\alpha\otimes\beta+\alpha\otimes\Delta\beta)  +(-\vu\otimes\vu + \F) \cdot *(\alpha\otimes\vN\beta))(t,\cdot )\\&  - \sum_{i,j} ( (h_{ij}) * ( \vN (\varphi \partial_i \partial_j G_d )*( \alpha\otimes\beta))) (t,\cdot) -  ( U_2 * (\alpha\otimes\beta)) (t,\cdot) .\end{split} \end{equation*} 

By Proposition \ref{poisson},  we conclude that $\tilde U *( \alpha\otimes\beta)(t,.)$ belongs to the space  $L^1(\Rd,(1+\vert x\vert^{-d-1}) $. Thus, it is a tempered distribution; as it is harmonic, it must be polynomial. The integrability in  $L^1(\Rd,(1+\vert x\vert^{-d-1}) $ implies that this polynomial is constant.

 If $\mathbb{F}$ belongs more precisely to  $L^1 ((0,T),L^1_{w_{d}}(\Rd))$
and $\vu$ belongs to $L^2 ((0,T),L^2_{w_{d}}(\Rd))$, we find that this polynomial belongs to  $L^1(\Rd, w_{d} \, dx)$, hence is equal to $0$.

Then, using the identity approximation $\Phi_\varepsilon = \frac 1{\varepsilon^4} \alpha	(\frac t\varepsilon)\beta(\frac x\varepsilon)$ and letting $\epsilon$ go to $0$,  we obtain a similar result for  $\tilde U$. Thus ${\bf S}=\vN p_\varphi +f(t)$, with $f(t)=0$ if   $\mathbb{F}$ belongs     to  $L^1 ((0,T),L^1_{w_{d}}(\Rd))$
and $\vu$ belongs to $L^2 ((0,T),L^2_{w_{d}}(\Rd))$.

As $f$ does not depend on $x$, we may take a function $\beta\in\mathcal{D}(\Rd)$ with $\int\beta\, dx=1$ and write $f=f*_x\beta$; we find that
$$ f(t)=\partial_t (\vu_0*\beta-\vu*\beta+\int_0^t \vu*\Delta\beta -(\vu\otimes\vu-\mathbb{F})\cdot*\vN\beta-p_\varphi* \vN\beta \, ds)=\partial_t g.$$
As $\partial_t\partial_j g=\partial_j f=0$ and $\partial_j g(0,.)=0$, we find that $g$ depends only on $t$; moreover, the formula giving $g$ proves that $g\in L^1((0,T))$.\Endproof\\

The proof of Theorem \ref{thbis} is classical and the result  is known as the extended Galilean invariance of the Navier--Stokes equations :\\

\pv{}  Let us suppose that 
\begin{equation*}
      \partial_t  \vu = \Delta \vu - ( \vu \cdot \vN ) \vu - \vN p_\varphi   - \frac{d}{dt} g(t),
\end{equation*}
  with  $g \in L^1((0,T))$. We define
  \begin{equation*}
      E(t)= \int_0^t g (\lambda) d\lambda \text{ and } \vw=\vu(t, x -E(t)) + g(t).
  \end{equation*}
We have  
  \begin{align*}
      \partial_t \vw =& \partial_t \vu(t, x-E(t))  - g(t) \cdot \vN \vu (t, x-E(t)) + \frac{d}{dt} g(t) \\
      =&\Delta \vu (t, x-E(t)) - [( \vu \cdot \vN ) \vu] (t, x-E(t)) - \vN p_\varphi (t, x-E(t)) -\frac{d}{dt} g(t)
      \\ &- g(t) \cdot \vN \vu (t, x-E(t))  + \frac{d}{dt} g(t) \\
      =&  \Delta \vw - ( \vw \cdot \vN ) \vw - \vN p_\varphi (t, x-E(t))  .
  \end{align*}

  If we define $q_\varphi(t,x)=p_\varphi(t,x-E(t))$, we find  that we have
  $$  q_\varphi=   {\displaystyle \sum_{i,j} } (\varphi \partial_i \partial_j G_d) * (w_iw_j)  +    {\displaystyle \sum_{i,j} }    \int (A_{i,j,\varphi}(x-y)-A_{i,j,\varphi}(-y))  (w_i(t,y)w_j(t,y) )\, dy.$$
  The theorem is proved. \Endproof

\section{Applications}
A consequence of Proposition \ref{poisson} is that we may define the Leray projection operator on the divergence of tensors that belong to $L^1((0,T),L^1(\Rd, w_{d+1}\, dx))$~: 
\begin{Definition} Let $\mathbb{H}\in  L^1((0,T),L^1(\Rd, w_{d+1}\, dx))$ and $\vw=\vN\cdot \mathbb{H}$. The Leray projection $\mathbb{P}(\vw)$ of $\vw$ on solenoidal vector fields is defined by
$$ \mathbb{P} \vw=\vw - \vN p_\varphi$$ where $\vN p_\varphi$ is the unique solution of $$-\Delta \vN p= \vN(\vN\cdot\vw)$$ such that
$$ \lim_{\tau\rightarrow +\infty} e^{\tau\Delta} \vN p=0.$$
\end{Definition}

A special form of the Navier--Stokes equations is then given by
$$(MNS) \quad  \partial_t\vu=\Delta\vu-\mathbb{P} \vN\cdot(\vu\otimes\vu-\mathbb{F}),\quad \vu(0,.)=\vu_0.$$ This leads to the integro-differential equation
$$ u=e^{t\Delta} \vu_0-\int_0^t e^{(t-s)\Delta}\mathbb{P} \vN\cdot(\vu\otimes\vu-\mathbb{F})\, ds.$$ The kernel of the convolution operator $e^{(t-s)\Delta}\mathbb{P} \vN\cdot$ is called the Oseen kernel; its study is the core of the method of mild solutions of Kato and Fujita \cite{FuK}. Thus, we will call equations (MNS)  a mild formulation of the Navier--Stokes equations.

The mild formulation together with the local Leray energy inequality has been as well a key tool for extending Leray's theory of weak solutions  in $L^2$ to the setting of weak solutions with infinite energy. We may propose a general definition of suitable Leray-type weak solutions :

\begin{Definition}[Suitable Leray-type solution]$\ $\\ Let $\mathbb{F}\in L^2((0,T), L^2(\Rd, \frac 1{(1+\vert x\vert)^{d+1}}))$ and $\vu_0\in L^2(\Rd, \frac 1{(1+\vert x\vert)^{d+1}})$ with $\vN\cdot \vu_0=0$. We consider the Navier--Stokes problem on $(0,T)\times\Rd$ :
  \begin{align*}
\partial_t\vu=& \Delta\vu-\mathbb{P}(\vu\otimes\vu-\mathbb{F}),\\ \vN\cdot\vu=0,& \quad\vu(0,.)=\vu_0.
  \end{align*}

A suitable Leray-type solution  $\vu$ of the Navier--Stokes equations is a vector field $\vu$ defined on $(0,T)\times\Rd$ such that :\begin{itemize}
\item[$\bullet$] $\vu$ is locally $L^2_t H^1_x$ on $(0,T)\times\Rd$
\item[$\bullet$] $\sup_{0<t<T} \int \vert\vu(t,x)\vert^2 \frac 1{(1+\vert x\vert)^{d+1}}\, dx<+\infty$
\item[$\bullet$] $\iint_{(0,T)\times\Rd} \vert\vN\otimes\vu (t,x)\vert^2 \frac 1{(1+\vert x\vert)^{d+1}}\, dx\, dt<+\infty$
\item[$\bullet$] the application $t\in [0,T)\mapsto \int \vu(t,x)\cdot \vw(x)\, dx$ is continuous for every smooth compactly supported vector field $\vw$
\item[$\bullet$] for every compact subset $K$ of $\Rd$, $\lim_{t\rightarrow 0} \int_K \vert \vu(t,x)-\vu_0(x)\vert^2\, dx=0$.
\item[$\bullet$] defining $p_\varphi$ as (the) solution of $-\Delta p_\varphi= \sum_{i,j}\partial_i\partial_j (u_iu_j-F_{i,j})$ given by Proposition \ref{poisson2},  $\vu$ is suitable in the sense of Caffarelli, Kohn and Nirenberg : there exists a non-negative locally bounded Borel measure $\mu$ on $(0,T)\times\Rd$ such that
$$ \partial_t(\frac{\vert\vu\vert^2}2)=\Delta (\frac{\vert\vu\vert^2}2)- \vert\vN\otimes\vu\vert^2-\vN\cdot ((\frac{\vert\vu\vert^2}2+p_\varphi)\vu)+ \vu\cdot(\vN\cdot\mathbb{F})
-\mu$$
\end{itemize}
\end{Definition}

\noindent{\bf Remarks : }\\ a) With those hypotheses, $p_\varphi$ belongs locally to $L^{3/2}_{t,x}$ and $\vu$ belongs locally to $L^3_{t,x}$ so that the distribution $(\frac{\vert\vu\vert^2}2+p_\varphi)\vu$ is well-defined.\\
b) Suitability is a local assumption. It has been introduced by Caffarelli, Kohn and Nirenberg in 1982 \cite{CKN} to get estimates on partial regularity for weak Leray solutions. If we consider a  solution of the Navier--Stokes equations on a small domain with no specifications on the behaviour of $\vu$ at the boundary, the estimates on the pressure (and the Leray projection operator) are no longer available. However, Wolf described in 2017 \cite{Wo17}  a local decomposition of the pressure into a term similar to the Leray projection of $\vN\cdot (\vu\otimes\vu)$ and a harmonic term; he could generalize the notion of suitability to this new description of the pressure. On the equivalence of various notions of suitability, see the paper by Chamorro, Lemari\'e-Rieusset and Mayoufi \cite{DLM}.\\
c) The relationship between the system (NS) and its mild formulation (MNS)  described in Theorem \ref{th} has been described by Furioli, Lemari\'e--Rieusset and Terraneo in 2000 \cite{FLT, LR02} in the context of uniformly locally square integrable solutions. See the paper by Dubois \cite{Dub}, as well.
\\

We list here a few examples to be found in the litterature :

\begin{enumerate}
\item Solutions in $L^2$ : in 1934, Leray \cite{Le34} studied the Navier--Stokes problem (NS) with an initial data $\vu_0\in L^2$ and a forcing tensor $\mathbb{F}\in L^2_tL^2_x$. He then obtained a solution $\vu\in L^\infty L^2\cap L^2\dot H^1$. Remark that this solution is automatically a solution of the mild formulation of the Navier--Stokes equations (MNS). Leray's construction by mollification provides suitable solutions.
\item Solutions in $L^2_{\rm uloc}$ : in 1999, Lemari\'e-Rieusset \cite{LR99, LR02}  studied the Navier--Stokes problem (MNS) with an initial data $\vu_0\in L^2_{\rm uloc}$ (and, later in \cite{LR16},  a forcing tensor $\mathbb{F}\in (L^2_tL^2_x)_{\rm uloc}$). He obtained (local in time) existence of a suitable solution $\vu$  on a small strip $(0,T_0)\times\Rd$ such that
$$ \sup_{x_0\in\Rd} \sup_{0<t<T_0} \int_{B(x_0,1)} \vert\vu(t,x)\vert^2\, dx <+\infty$$ and
$$ \sup_{x_0\in\Rd} \int_0^{T_0}\int_{B(x_0,1)} \vert\vN\otimes\vu(t,x)\vert^2\, dx <+\infty. $$ Remark that we have $\vu\in L^2((0,T_0), L^2(\Rd,\frac 1{(1+\vert x\vert)^{d+1}}\, dx))$ but $\vu$ does not belong to  $ L^2((0,T_0), L^2(\Rd,\frac 1{(1+\vert x\vert)^{d}}\, dx))$; thus, in this setting, problems (NS) and (MNS) are not equivalent.

Various reformulations of local Leray solutions in $L^2_{\rm uloc}$ have been provided, such as  Kikuchi and Seregin in 2007 \cite{KS07}  or Bradshaw and Tsai in 2019 \cite{BT19}. The formulas proposed for the pressure, however, are actually equivalent, as they all imply that $\vu$ is solution to the (MNS) problem.

In the case of dimension $d=2$, Basson \cite{Ba06}  proved in 2006 that the solution $\vu$ is indeed global (i.e. $T_0=T$) and that, moreover, the solution is unique.
\item Solutions in a weighted Lebesgue space : in 2019, Fern\' andez-Dalgo and Lemari\'e--Rieusset  \cite{PF_PG} considerered data $\vu_0\in L^2(\Rt,  w_\gamma\, dx)$ and $\mathbb{F}\in L^2((0,+\infty), L^2(\Rt, w_\gamma\, dx))$ with $0<\gamma\leq 2$. They proved  (global in time) existence of a suitable solution $\vu$    such that, for all $T_0<+\infty$, 
$$  \sup_{0<t<T_0} \int \vert\vu(t,x)\vert^2 \, w_\gamma(x)\, dx <+\infty$$ and
$$  \int_0^{T_0}\int \vert\vN\otimes\vu(t,x)\vert^2\, w_\gamma(x)\\, dx <+\infty. $$ [Of course, for such solutions, (NS) and (MNS) are equivalent.]  They showed that, for $\frac 4 3<\gamma\leq 2$, this frame of work is well adapted to the study of discretely self-similar solutions with locally $L^2$ initial value, providing a new proof of the results of Chae and Wolf in 2018 \cite{CW18} and of Bradshaw and Tsai in 2019 \cite{BT19a}.
\item Homogeneous Statistical Solutions : in 1977, Vishik and Fursikov \cite{VF} considered the (MNS) problem with a random initial value $\vu_0(\omega)$. The statistics of the initial distributions were supposed to be invariant though translation of the arguments of $\vu_0$ : for every Borel subset $B$ of $L^2_{\rm loc}(\Rt)$ and every $x_0\in\Rt$,
$$ Pr( \vu_0(\cdot - x_0)\in A)=Pr(\vu_0\in A).$$
Another assumption was that $\vu_0$ has a bounded mean   energy density  :
$$ e_0=\mathbb{E}\left(\frac{\int_{\vert x\vert\leq 1} \vert \vu_0\vert^2\, dx}{\int_{\vert x\vert\leq 1}\, dx}\right)<+\infty.$$  Then $$ Pr( \vu_0\in L^2\text{ and } u\neq 0)=0$$ while, for any $\epsilon>0$, 
$$ Pr( \int \vert \vu_0\vert^2\frac 1{(1+\vert x\vert)^{3+\epsilon}}\, dx<+\infty)=1.$$ In \cite{VF88}, they  constructed a solution $\vu(t,x,\omega)$ that solved the Navier--Stokes equation for almost every initial value $\vu_0(\omega)$, and the solution belonged almost surely to $L^\infty_t L^2_x(\frac 1{(1+\vert x\vert)^{3+\epsilon}}\, dx )$ with $\vN\otimes\vu\in L^2_tL^2_x(\frac 1{(1+\vert x\vert)^{3+\epsilon}}\, dx )$. 

In 2006, Basson \cite{Ba06a} gave a precise description of the pressure in those equations (which is equivalent to our description through the Leray projection operator) and proved the suitability of the solutions.
\end{enumerate}

\section{The space $B_\gamma^2$.}

Instead of dealing with weighted Lebesgue spaces, one may deal with a kind of local Morrey space, the space $B_\gamma^2$.

\begin{Definition}
For $\gamma\geq 0$, define $w_\gamma(x)=\frac 1{(1+\vert x\vert)^\gamma}$ and $L^p_{w\gamma}=L^p(\Rd, w_\gamma(x)\, dx)$.

For $1\leq p < +\infty $, we denote $B^p_{\gamma}$ the Banach space  of all functions $ u \in L^p_{ \rm loc} $ such that :
\begin{equation*}
    \| u  \|_{B_{\gamma}^p} = \sup_{R \geq 1}  ( \frac{1}{R^{\gamma }} \int_{B(0,R)} |u|^{p} \, dx  )^{1/p} < +\infty.
\end{equation*}
 
Similarly,  $B^p_\gamma L^p (0,T)$ is  the Banach space of all functions $u \subset (L^p_t L^p_x)_{ \rm loc}$ such that
\begin{equation*}
    \| u \|_{B^p_\gamma L^p (0,T)} = \sup  \left( \frac{1}{R^\gamma} \int_0^T \int |u|^p \right)^{\frac{1}{p}} dx \, ds.
\end{equation*}
\end{Definition}

\begin{Lemma}
\label{li}
Let $\gamma \geq 0$ and $  \gamma < \delta < + \infty $, we have the continuous embedding  $L^p_{w_{\gamma}} \hookrightarrow B^p_{\gamma,0 } \hookrightarrow B^p_{\gamma } \hookrightarrow L^p_{w_{\delta}} $, where $B^p_{\gamma,0 } \subset B^p_{\gamma }$ is the subspace of all functions $u \in B^p_{\gamma }$ such that $  \lim_{ R \to +\infty }  \frac{1}{ R^{\gamma }} \int_{B(0,R)} | u (x) |^p \, dx =0 $. 
\end{Lemma}

\pv
Let $u \in L^p_{w_\gamma}$. We verify easily that $ \| u \|_{B^p_{\gamma}} \leq 2^{\gamma/p} \| u \|_{ L^p_{w_{\gamma}} }$ and we see that
\begin{equation*}
       \frac{1}{ R^{\gamma }} \int_{|x|\leq R} |u|^{p} \, dx  = \int_{ |x| \leq R } \frac{|u|^{p}}{(1+|x|)^{\gamma}}   \frac{(1+ |x| )^{\gamma}}{R^{\gamma }} \, dx 
\end{equation*}
converges to zero when $    R \to + \infty $ by dominated convergence, so $L^p_{w_{\gamma}} \hookrightarrow B^p_{\gamma,0 }$. To demonstrate the other part, we estimate
\begin{align*}
       \int \frac{ |u|^{p}}{  (1+|x|)^{\delta }} \, dx &= \int_{|x|\leq 1} \frac{ |u|^{p}}{  (1+|x|)^{\delta }} \, dx + \sum_{n \in \mathbb{N}} \int_{2^{n-1} \leq |x| \leq 2^{n}} \frac{ |u|^{p}}{  (1+|x|)^{\delta }} \, dx \\
       & \leq \int_{|x|\leq 1} |u|^{p} \, dx + \sum_{n \in \mathbb{N}} \frac{1}{(1+2^{n-1})^{\delta}}  \int_{2^{n-1} \leq |x| \leq 2^{n}}  |u|^{p} \, dx  \\
       & \leq \int_{|x|\leq 1} |u|^{p} \, dx + c \sum_{n \in \mathbb{N}}  \frac{1}{2^{\delta n } }   \int_{2^{n-1} \leq |x| \leq 2^{n}}  |u|^{p} \, dx \\
       &  \leq  ( 1 + c \sum_{n \in \mathbb{N}}  \frac{1}{2^{(\delta -\gamma ) n } }   ) \sup_{R \geq 1}   \frac{1}{R^\gamma } \int_{|x| \leq R } |u|^{p} \, dx  ,
\end{align*}
thus, $ B^p_{\gamma} \subset L^p_{w_{\delta}} $.
 \Endproof\\
 
 \noindent{\bf Remark : } Similarly,  for all $\delta > \gamma$, $ B^p_\gamma L^p (0,T) \subset L^p ((0, T), L^p_{w_\delta})$.

\begin{Proposition}
\label{ti}
The space $B^p_\gamma$ can be obtained by interpolation, 
\begin{equation*}
    B^p_\gamma = [L^p, L^p_{w_{\delta}} ]_{\frac{\gamma}{\delta}, \infty}
\end{equation*}
for all $0 < \gamma < \delta < \infty$,
and the norms
\begin{equation*}
    \| \cdot \|_{B^p_\gamma} \phantom{space} \text{and} \phantom{space} \| \cdot \|_{[L^p, L^p_{w_{\delta}} ]_{\frac{\gamma}{\delta}, \infty}}
\end{equation*}
are equivalent.
\end{Proposition}

\pv
Let $f \in B^p_\gamma$. For $A<1$, we write $f_0 = 0$ and $f_1 = f$, then we have $f=f_0+f_1$ and $\| f_1 \|_{L^p_{w_\delta}} \leq C A^{\frac{\gamma}{\delta}-1} \| f \|_{B^p_{\gamma}} $. 

For $A>1$, we let $R=A^{\frac{p}{\delta}} > 1$. We write $f_0 = f \mathds{1}_{|x| \leq R}$ and $f_1 = f \mathds{1}_{|x| > R}$, then
\begin{equation*}
    \| f_0 \|_p \leq C \| f \|_{B^p_\gamma} R^{\frac{\gamma}{p}} = C A^{\frac{\gamma}{\delta}} \| f \|_{B^p_{\gamma}}
\end{equation*}
and
\begin{align*}
    \| f_1 \|^p_p &= \sum_{n \in \mathbb{N}} \int_{2^{n-1}R \leq |x| \leq 2^{n} R} \frac{ |u|^{p}}{  (1+|x|)^{\delta }} \, dx \\
    & \leq C R^{\gamma - \delta } \sum_{n \in \mathbb{N}} \frac{1}{2^{(\delta -\gamma)j}} \| f \|^p_{B^p_\gamma} \\
    &= C A^{(\frac{\gamma}{\delta}-1)p} \| f \|^p_{B^p_{\gamma}}
\end{align*}
Thus, $B^p_\gamma \hookrightarrow [L^p, L^p_{w_{\delta}} ]_{\frac{\gamma}{\delta}, \infty}$.

Let $f \in [L^p, L^p_{w_{\delta}} ]_{\frac{\gamma}{\delta}, \infty}$, then there exist $c>0$ such that for all $A>0$, there exist $f_0 \in L^p$ and $f_1 \in L^p_{w_{\delta}}$ so that $f=f_0 + f_1$,
\begin{equation*}
     \| f_0 \|_p \leq  c A^{\frac{\gamma}{\delta}} 
\phantom{space}
and
\phantom{space}
    \| f_1 \|_{L^p_{w_\delta}} \leq
     c A^{\frac{\gamma}{\delta}-1} .
\end{equation*}
For $j \in \mathbb{N}$ we take $A= 2^{\frac{j \delta}{p}}$, then
\begin{align*}
    & \frac{1}{2^{j \gamma }} \int_{|x| < 2^j} |f|^p \, dx \\
    & \leq C \left(   \frac{1}{2^{j \gamma }} \int_{|x| < 2^j} |f_0|^p \, dx + \frac{1}{2^{j \gamma }} \int_{|x| < 2^j} |f_1|^p \, dx \right) \\
    & \leq C \left(    \frac{1}{2^{j \gamma }}  \|f_0\|_p^p+ \frac{C}{2^{j \gamma }} \int_{|x| \leq 1}   \frac{|f_1|^p}{(1+|x|)^{ \delta}} \, dx + C \sum_{k=1}^j \frac{2^{k \delta}}{2^{j \gamma }}  \int_{2^{k-1} < |x| < 2^k} \frac{|f_1|^p}{(1+|x|)^{ \delta}} \, dx  \right) \\
    & \leq C \left(    \frac{1}{2^{j \gamma }}  \|f_0\|_p^p+ C' 2^
    {j(\delta -\gamma)} \| f_1 \|^p_{L^p_{w_\delta}}  \right)  \\
    & \leq C''  
\end{align*}
which implies $\sup_{j \in \mathbb{N}} \frac{1}{2^{j \gamma }} \int_{|x| < 2^j} |f|^p \, dx < +\infty $, so $\sup_{R\geq 1 } \frac{1}{R^\gamma} \int_{|x| < R} |f|^p \, dx < +\infty $. \Endproof\\

Thus, we can see that the local Morrey spaces $B^p_\gamma$ are very close to the weighted Lebesgue spaces $L^p_{w_\gamma}$. Indeed, the methods and results of Fern\' andez-Dalgo and Lemari\'e--Rieusset  \cite{PF_PG}  can be easily extended to the setting of local Morrey spaces in dimension $d=2$ or $d=3$  : considering data $\vu_0\in B^2_\gamma(\Rd)$ and $\mathbb{F}\in (B^2_\gamma L^2)(0,T)(\Rd)$ with $0<\gamma\leq 2$, one gets  (local in time) existence of a suitable solution $\vu$ for the (MNS) system  on a small strip $(0,T_0)\times\Rd$  such that $\vu\in L^\infty((0,T_0), B^2_\gamma)$ and $\vN\otimes\vu \in (B^2_\gamma L^2)(0,T_0)$. 

The case of $\gamma=2$ deserves some comments. In the case $d=3$, the results is slightly more general than the results in \cite{PF_PG}, as the class $B^2_2$ is larger than the space $L^2_{w_2}$. Equations in $B^2_2$ have been very recently discussed by Bradshaw, Kukavica and Tsai  \cite{BTK}. The case $d=2$ is more intricate. Indeed, while the Leray projection operator is bounded on $B^2_2(\Rt)$ (by interpolation with $L^2$ and $L^2_{w_\delta}$ with $2<\delta<3$, the Riesz transforms being bounded on 
 $L^2_{w_\delta}$ by the theory of Muckenhoupt weights), this is no longer the fact on $B^2_2(\Rdo)$. Thus, one must be careful in the handling of the pressure. This has been done by Basson in his Ph. D. thesis in 2006 \cite{Ba06}.
 
 Local Morrey spaces $B^2_d$ occur naturally in the setting of homogeneous statistical solutions. By using an ergodicity argument, Dostoglou \cite{Do01} proved in 2001 that, under the assumptions of Vishik and Fursikov \cite{VF}, we have
 $$ Pr( \vu_0(.,\omega)\in B^2_d(\Rd))=1.$$ Thus, the solutions of Vishik and Fursikov live in a smaller space than $L^2_{w_{d+\epsilon}}$.


\begin{thebibliography}{HD}

 \bibitem{Ba06}	  A. Basson, \emph{Solutions spatialement homog\`enes adapt\'ees des \'equations de Navier--Stokes}, Th\`ese, Universit\'e d'\'Evry, 2006.

 \bibitem{Ba06a}	  A. Basson, \emph{Homogeneous Statistical Solutions and Local Energy
Inequality for 3D Navier--Stokes Equations}, Commun. Math. Phys. 266 (2006), 17--35.
 
 
\bibitem{BT19a}	Z. Bradshaw and T.P. Tsai, \emph{Discretely self-similar solutions to the Navier-Stokes
equations with data in $L^2_{\rm loc}$}, to appear in  Analysis and PDE.	

\bibitem{BT19} Z. Bradshaw and T.P. Tsai, \emph{Global existence, regularity, and uniqueness of infinite energy
solutions to the Navier--Stokes equations}, arXiv:1907.00256.

\bibitem{BTK}	Z. Bradshaw, I. Kukavica and T.P. Tsai, \emph{Existence of global weak solutions to the Navier-Stokes equations in weighted spaces}, arXiv:1910.06929v1

\bibitem{CKN} L. Caffarelli, R. Kohn and  L. Nirenberg, \emph{Partial regularity of suitable weak solutions of the Navier--Stokes
equations},  Comm. Pure Appl. Math., 35 (1982), 771--831.

\bibitem{CW18} 	D. Chae and J. Wolf,  \emph{Existence of discretely self-similar solutions to the Navier-Stokes
equations for initial value   in $L^2_{\rm loc}(\mathbb{R}^3)$}, Ann. Inst. H. Poincar\'e  Anal. Non Lin\'eaire  {35} (2018), 1019--1039.

\bibitem{DLM} D. Chamorro, P.G.  Lemari\'e--Rieusset, and K.  Mayoufi, \emph{The role
of the pressure in the partial regularity theory for weak solutions of the Navier--Stokes
equations},  Arch. Rat. Mech. Anal. 228.1 (2018), 237--277.

\bibitem{PF_PG} P.G Fern\'andez-Dalgo, P.G.  Lemari\'e--Rieusset, \emph{Weak solutions for Navier--Stokes equations with
initial data in   weighted $L^2$ spaces.}, preprint, 2019.

\bibitem{Do01} S. Dostoglou, \emph{Homogeneous measures and spatial ergodicity of the Navier--Stokes equations}, preprint, 2002.

\bibitem{Dub} S. Dubois, \emph{What is a solution to the Navier--Stokes equations?}, 
C. R. Acad. Sci. Paris, Ser. I 335 (2002,) 27--32.
\bibitem{FuK}
  H. Fujita and T. Kato, \emph{On the non-stationary Navier-Stokes system}, Rendiconti Seminario Math. Univ. Padova {32} (1962),  243–260.  
  
  	\bibitem{FLT} 
	G. Furioli, P.G. Lemari\'e-Rieusset and E. Terraneo. \emph{Unicit\'e dans
$\text{L}^3(\Rt)$ et d'autres espaces limites pour
Navier--Stokes}, Revista Mat.
Iberoamericana  16 (2000), 605--667.


\bibitem{KS07} N. Kikuchi and G. Seregin, \emph{Weak solutions to the Cauchy problem for the Navier--Stokes equations satisfying the local energy inequality}. Nonlinear equations and spectral
theory, 141--164, Amer. Math. Soc. Transl. Ser. 2, 220, Amer. Math. Soc., Providence,
RI, 2007.



   \bibitem{LR99}
P.G.  Lemari\'e--Rieusset, \emph{Solutions faibles d'\'energie infinie pour les   \'equations de {N}avier--Stokes dans $\mathbb{R}^{3}$}, C. R. Acad. Sci. Paris, Serie I. {328} (1999), 1133--1138.
	
    \bibitem{LR02}  P.G. Lemari\'e-Rieusset, \emph{Recent developments in the Navier--Stokes problem}, CRC Press, 2002.



    \bibitem{LR16}
P.G.  Lemari\'e--Rieusset, \emph{The Navier--Stokes problem in the 21st century}, Chapman \& Hall/CRC, (2016).

    \bibitem{Le34}
J. Leray, \emph{Essai sur le mouvement d'un fluide
visqueux emplissant l'espace}, Acta Math. {63} (1934), 193--248.

\bibitem{VF}   M.I.  Vishik and A. V.  Fursikov, \emph{Solutions statistiques homog\`enes des syst\`emes diff\'erentiels paraboliques
et du syst\`eme de Navier-Stokes}. Ann. Scuola Norm. Sup. Pisa, série IV, IV (1977), 531--576.

\bibitem{VF88}   M.I.  Vishik and A. V.  Fursikov, \emph{Mathematical Problems of Statistical Hydromechanics}, Dordrecht:
Kluwer Academic Publishers, 1988.

\bibitem{Wo17} J. Wolf, \emph{On the local pressure of the Navier--Stokes equations and related systems},     Adv. Differential Equations
 22 (2017), 305--338.
\end{thebibliography}
\end{document}